\documentclass[12pt,utf]{extarticle}

\usepackage{caption}
\usepackage{enumerate}
\usepackage{t1enc}
\usepackage{babel}
\usepackage{amsthm}
\usepackage{amsmath}
\usepackage{amsfonts}
\usepackage{mathrsfs}
\usepackage[cp1250]{inputenc}
\usepackage{amsthm}
\usepackage{amsmath}
\usepackage{amssymb}
\usepackage{amsfonts}
\usepackage{indentfirst}
\usepackage{textcomp}
\usepackage{mathabx}
\usepackage{bbm}
\usepackage{tikz}
\usepackage{xcolor}
\usepackage{extsizes}
\usepackage[textwidth=18.5cm,textheight=23cm,left=3.5cm,voffset=-1cm,,footskip=1.5cm,headheight=0cm]{geometry}
\usepackage{fancyhdr}
\usepackage{yfonts}
\usepackage{float}
\usepackage{pgf}
\usepackage{pgfplots}

\usepackage{mathtools}

\pgfplotsset{compat=1.12}

\title{\LARGE \bf Markov control of continuous time Markov processes with  long run functionals by time discretization*
}

\author{Lukasz Stettner$^{1}$ 
\thanks{*This work was supported by NCN grant 2024/53/B/ST1/00703}
\thanks{$^{1}$Lukasz Stettner is with the Institute of Mathematics Polish Academy of Sciences,
        Sniadeckich 8, 00-656 Warsaw, Poland,
        {\tt\small stettner@impan.pl}}%
}

\begin{document}

\maketitle
\thispagestyle{empty}
\pagestyle{empty}

\begin{abstract}

In the paper we study continuous time controlled  Markov processes using discrete time controlled Markov processes. We consider long run functionals: average reward per unit time or long run risk sensitive functional. We also investigate stability of continuous time functionals with respect to pointwise convergence of Markov controls.  
\end{abstract}

\section{INTRODUCTION}
Assume that state space $E$ is Polish with Borel $\sigma$-field ${\cal E}$, although in particular examples we shall consider ${\cal E}=R^d$ or a bounded convex subset of $R^d$. We have also a compact set of control parameters $U$ and a family ${\cal U}$ of Borel measurable mappings $u: E \mapsto U$ called later Markov controls. On a probability space $(\Omega,F,(F_t),P)$,  for each $u\in {\cal U}$ we are given a continuous time controlled Markov process $(X_t^u)$ with transition operator $P^{u}_t(x,dy)$ for $x\in E$ and control $u(X_t^u)$ at generic time $t$. We consider a natural pointwise convergence topology on ${\cal U}$, which means that $u_n\in {\cal U}$ converges to $u\in {\cal U}$ whenever $u_n(x)\to u(x)$ as $n\to \infty$ for each $x\in E$.
Then we consider discrete time approximations $(X_t^{(h),u})$ of $(X_t^u)$ which is a discrete Markov process $X_{nh}^{(h),u}$ at generic moments $nh$ such that 
$X_t^{(h),u}=X_{\left[{t\over h}\right]h}^{(h),u}$, where $\left[{t\over h}\right]$ is the integer part of ${t \over h}$   and $X_{nh}^{(h),u}$ has transition operator $P^{(h),u(X_{nh}^{(h),u})}(X_{nh}^{(h),u},\cdot)$. This means that while process $(X_t^u)$ is controlled at each time $t$ using $u(X_t^u)$, its discrete time approximation $X_{nh}^{(h),u}$ is controlled at moments $nh$ using $u(X_{nh}^{(h),u})$. To be more precise consider our main example.

{\bf Example 1.} Assume for $u\in {\cal U}$ we have the following equation in $R^d$ 
\begin{equation}\label{equat1}
X_t^u=x_0 + \int_0^t b(X_s^u,u(X_s^u))ds + \int_0^t \sigma(X_s^u)dW_s,
\end{equation}
where $(W_t)$ is a Brownian motion, $|b(x,a)-b(y,a)|+\|\sigma(x)-\sigma(y)\|\leq K_R |x-y|$
for $a\in U$, $|x|, |y| \leq R$, $|b(x,a)|^2+\|\sigma(x)\|^2\leq K(1+|x|^2)$ and $\xi^T \sigma(x)\sigma^T(x)\xi \geq {1\over K_R}|\xi|^2$
for $\xi\in R^d$, $|x|\leq R$ and any $R>0$. By Theorem 2.2.12 of \cite{An} for each $u\in {\cal U}$ there is a unique strong solution to the equation (\ref{equat1}).    
Our discrete approximation with discretization step $h$ is defined as 
\begin{equation}\label{equat2}
X_{(n+1)h}^{(h),u}= X_{nh}^{(h),u}+\int_{nh}^{(n+1)h} b(X_s^{(h),u},u(X_{nh}^{(h),u}))ds +  \int_{nh}^{(n+1)h} \sigma(X_s^{(h),u})dW_s 
 \end{equation}
for $n=1,2 \ldots$ and $X_0^{(h),u}=x$. Since we have a unique strong solution on each time interval $[nh,(n+1)h]$ we have well defined process $(X_t^{(h),u})$. In what follows we shall consider a general case introducing a number of assumptions which are mainly satisfied by the model considered in this example. 

In the paper we want to maximize the following functionals: {\it average reward per unit time} 

\begin{equation}\label{ACc}
J_x(u)=\liminf_{t\to \infty}{1\over t} E_x^u\left\{\int_0^t c(X_s^u,u(X_s^u))ds\right\},
\end{equation}
 for a bounded measurable function $c: E\times U\mapsto R$, continuous with respect to the second (control) parameter,
  
and its discrete time approximation
\begin{equation}\label{ACd}
J_x^h(u)=\liminf_{n\to \infty}{1\over n h} E_x^u\left\{\sum_{i=0}^{n-1} h c(X_{ih}^{(h),u},u(X_{ih}^{(h),u}))\right\},
\end{equation}

\noindent 
{\it long run risk sensitive} with risk factor $\alpha<0$

\begin{equation}\label{RACc}
 I_x^\alpha(u)=\liminf_{t\to \infty} {1\over \alpha} {1\over t} \ln E_x^u\left\{e^{\alpha \int_0^t c(X_s^u,u(X_s^u))ds}\right\},
\end{equation}

and its discrete time approximation
\begin{equation}
 I_x^{\alpha,h}(u)=\liminf_{n\to \infty} {1 \over \alpha} {1\over n h} \ln E_x^u\left\{e^{\alpha h \sum_{i=0}^{n-1} c(X_{ih}^{(h),u},u(X_{ih}^{(h),u}))}\right\}.
\end{equation} 
Risk sensitive functionals are important since they measure not only expected value of the reward but also other moments of the reward including variance with weight $\alpha$, which his considered as a measure of risk (see \cite{PS}, \cite{LS3}, \cite{LS4}). 
We want to show that under suitable assumptions $J_x^h(u)\to J_x(u)$ and 
$I_x^{\alpha,h}(u)\to I_x^{\alpha}(u)$ as $h\to 0$. Then we consider stability of continuous time functionals i.e. we using discrete approximation  show that whenever $u_n\to u$ then also $J_x(u_n)\to J_x(u)$ and  $I_x^\alpha(u_n)\to  I_x^\alpha(u)$ as $n\to \infty$. 

The paper generalizes and extends \cite{LS0}, where only discrete time was considered. Usually we have a continuous time model which we control using discrete time inputs. In the paper we want to justify such procedure. Practically we use piecewise constant controls in discrete time moments, which we expect to be good, feasible approximation of real world model. Notice that such models can not be approximated using weak convergence technics considered in \cite{KD}.  Average reward per unit time problem is considered in full generality considering Lyapunov function $V$, which allows us to have unbounded reward function $c$ and consequently we obtain a number of results in norms weighted by $V$. The studies of long run risk sensitive functionals are practically restricted to compact state spaces for which we consider nondegenerate diffusions, possibly with jumps, in regular bounded sets.    

\section{Average reward per unit time problem}

We shall need the following assumption: 

\medskip

\noindent
(ER) {\it for each $u\in {\cal U}$ process $(X_t^u)$ is aperiodic and ergodic in the sense that it has a unique invariant measure $\mu^u$.}

\medskip

In what follows we shall consider discrete time approximations with $h=2^{-m}$, and to simplify notations we shall denote process $X_{nh}^{(h),u}$ by $X_{n2^{-m}}^{(m),u}$. We assume that 

\medskip

\noindent
(ERd) {\it for each $m\in N$ and the process $(X_{n2^{-m}}^{(m),u})$ is aperiodic and ergodic.}

\medskip
  
Furthermore we assume that 

\medskip

\noindent
(UEd) {\it for each $u\in {\cal U}$ there is $\rho\in (0,1)$ and function $V:E\to [1,\infty)$ such that for $x,x'\in E$ and $m \in N$
\begin{equation}\label{ergcond}
\int_E V(y)|P_1^{(m),u}(x,dy)-P_1^{(m),u}(x',dy)|\leq \rho \left[V(x)+V(x')\right].
\end{equation}}

\medskip

Above introduced $V$ is called sometimes a {\it Lyapunov function}. Using $V$ we consider the norm $\|f\|_V:=\sup_{x\in E} {|f(x)|\over V(x)}$ for Borel measurable functions $f$ and define the space $B_V$ as the space of Borel measurable functions $f$ with finite norm  $\|f\|_V$. Similarly in the space of finite signed measures $M(E)$ we consider the norm 
$\|\nu\|_V:= sup_{f\in B_V, \|f\|_V\leq 1} |\int _E f(x)\nu(dx)|$.

The condition (\ref{ergcond}) was introduced by Kartashov (see \cite{Ka} and also \cite{HLL}) and has the following important consequences

\noindent
{\bf Lemma 1.} {\it  If there is $x^*\in E$ such that $P_{n}^{(m),u}V(x^*)<\infty$ then under (UEd) there is a unique invariant measure $\mu_m^u$ for the Markov process $(X_{n}^{(m),u})$ and
\begin{equation}\label{impin}
{\|P_n^{(m),u}(x,\cdot)-\mu_m^u(\cdot)\|_V \over V(x)}\leq \rho^n[1+ {P_1^{(m),u}V(x^*)+\rho V(x^*) \over 1-\rho}].
\end{equation}}

\noindent{\bf Proof.} It follows from Theorem 7.3.14 of \cite{HLL}. 

\rightline{$\Box$}

Assume  

\medskip

\noindent
(FPV) {\it we have $\sup_m \sup_{x\in E}{P_{n}^{(m),u}V(x)\over V(x)}<\infty$ for each $x\in E$}

\medskip

We immediately have

\noindent
{\bf Corollary 1.} {\it If $\sup_m P_1^{(m),u}V(x^*)<\infty$ for some $x^*\in E$ then the bound in (\ref{impin}) is uniform with respect to $m$ and consequently we have (FPV). Assuming furthermore (ERd) we have that $\mu_m(\cdot)$ is a unique invariant measure for the process $(X_{n2^{-m}}^{(m),u})$.}

Denote by $P(E)$ the set of probability measures on $E$ and let $P_V(E):=\left\{\nu\in P(E): \|\nu\|_V<\infty \right\}$. In what follows we shall need the following technical Lemma

\noindent
{\bf Lemma 2.}\label{techlem1} 
{\it Assume that  for $\nu_n, \nu\in P_V(E)$ we have $\|\nu_n - \nu\|_V\to 0$ and for $f_n, f\in B_V$ with  $\|f_n\|_V$ bounded we have
${f_n(x)-f(x)} \to 0$ for each $x\in E$.  Then $\nu_n(f_n)\to \nu(f)$.}
  
\noindent{\bf Proof.}
Without loss of generality we may assume that $\|f_n\|_V\leq 1$. Then also $\|f\|_V\leq 1$ and we have
\begin{eqnarray}
&&|\nu_n(f_n)-\nu(f)|\leq |\nu_n(f_n)-\nu(f_n)|+|\nu(f_n)-\nu(f)\|\leq \nonumber \\
&& \|\nu_n-\nu\|_V + |\nu((g_n-g)V)|\to 0  
\end{eqnarray}
as $n\to \infty$, with $g_n={f_n\over V}$, $g={f\over V}$ and where the last convergence follows from the dominated convergence theorem. 

\rightline{$\Box$}

Assume 

\medskip

\noindent
(Conv) {\it for each $u\in {\cal U}$ and $x\in E$ we have $\|P_1^{(m),u}(x,\cdot)-P_1^{u}(x,\cdot)\|_V\to 0$ as $m\to \infty$. }

\medskip

We have

\noindent
{\bf Proposition 1.} {\it Under (Conv) and (FPV) for each $n\in N$ and $x\in E$ we have
\begin{equation}\label{convn}
\|P_n^{(m),u}(x,\cdot)-P_n^{u}(x,\cdot)\|_V\to 0
\end{equation} as $m\to \infty$.} 
   
\noindent{\bf Proof.} We use induction. For $n=1$ (\ref{convn}) is satisfied by (Conv). Assume that we have (\ref{convn}) for $n$. Then by (FPV) we have that there is $K\geq 0$ such that 
\begin{equation}
\sup_{f\in B_V, \|f\|_V\leq 1} \sup_{x\in E}{|P_{n}^{(m),u}(x,f)|\over V(x)}= \sup_{x\in E}{P_{n}^{(m),u}(x,V)\over V(x)}\leq K<\infty
\end{equation}
and therefore 
\begin{eqnarray}
&&\sup_{f\in B_V, \|f\|_V\leq 1} |P_{n+1}^{(m),u}(x,f)-P_{n+1}^{u}(x,f)|\leq \sup_{f\in B_V, \|f\|_V\leq 1} \left[|\int_E P_{n}^{(m),u}(y,f)(P_{1}^{(m),u}(x,dy)- \right. \nonumber \\
&&P_{1}^{u}(x,dy))| +|\int_E (P_{n}^{(m),u}(y,f)- \left. P_{n}^{u}(y,f))P_{1}^{u}(x,dy)|\right] \leq  K \|P_n^{(m),u}(x,\cdot)- \nonumber \\
&& P_n^{u}(x,\cdot)\|_V +\int_E \|P_n^{(m),u}(y,\cdot)-P_n^{u}(y,\cdot)\|_V P_{1}^{u}(x,dy)\nonumber 
\end{eqnarray}
and by induction hypothesis and dominated convergence we have that (\ref{convn}) for $n+1$ follows. 

\rightline{$\Box$}

Using Proposition 1 to (UEd) and then Lemma 1 we immediately obtain

\noindent
{\bf Corollary 2.} {\it Under (UEd), (FPV), (ER) and (Conv) we have
\begin{equation}
\int_E V(y)|P_1^{u}(x,dy)-P_1^{u}(x',dy)|\leq \rho \left[V(x)+V(x')\right]
\end{equation} 
and 
\begin{equation}
\|P_n^{u}(x,\cdot)-\mu^u(\cdot)\|_V \leq \rho^n V(x)[1+ {P_1^{u}V(x^*)+\rho V(x^*) \over 1-\rho}]
\end{equation}
where $\mu^u$ is a unique invariant measure for $(X_t^u)$.}

We can now rewrite the functional (\ref{ACd}) with $h=2^{-m}$ in the form 

\begin{equation}\label{ACdp}
J_x^m(u)=\liminf_{n\to \infty}{1\over n } E_x^u\left\{\sum_{i=0}^{n-1} C_m(X_{i}^{(m),u},u))\right\},
\end{equation}
where 
\begin{equation}
C_m(x,u):=E_x^u\left\{\sum_{i=0}^{2^{m}-1}2^{-m}c(X_{i 2^{-m}}^{(m),u},u(X_{i 2^{-m}}^{(m),u}))\right\}.
\end{equation}
 with a continuous time analog
\begin{equation}
C(x,u):= E_x^u\left\{\int_0^1 c(X_s^u,u(X_s^u))ds\right\} 
\end{equation}
We shall assume that

\medskip

\noindent
(CCon) {\it $C_m, C\in B_V$ and for each $x\in E$ we have that $\|C_m\|_V$ is bounded and $|C_m(x,u)-C(x,u)|\to 0$ for $x\in E$ and $u\in {\cal U}$, as $m\to \infty$.}  

\medskip

Notice that in this section we allow $c$ to be unbounded, we require only that $c\in B_V$ as in (CCon).
We have

\noindent
{\bf Theorem 1.} {\it Under (Conv), (FPV), (CCon), (ER) and (ERd) we have that 
\begin{equation}\label{invconv}
\|\mu_m^u-\mu^u\|_V\to 0,
\end{equation}
\begin{eqnarray}
&&J_x^m(u)=\int_E C_m(x,u) \mu_m^u(dx)= \int_E c(x,u(x))\mu_m^u(dx) \to \int_E C(x,u)\mu^u(dx)= \nonumber \\
&&\int_E c(x,u(x)\mu^u(dx)=J_x(u)
\end{eqnarray}
as $m\to \infty$.}

\noindent{\bf Proof.} (\ref{invconv}) follows from (\ref{impin}), (\ref{convn}) and Corollary 2. By 
Lemma 2 and (Ccon) we have that $\mu_m^u(C_m)\to \mu^u(C)$. Now from (ERd) we have that $\mu_m^u(C_m)=\int_E c(x,u(x)\mu_m^u(dx)$, while from (ER) we have that $\mu^u(C)=\int_E c(x,u(x)\mu^u(dx)$, which completes the proof. 

\rightline{$\Box$}
  
To study continuity of the cost functional $J_x^h(u)$ with respect to $u\in {\cal U}$ we shall need the following assumption

\medskip

\noindent
(uCont) {when $u_n\to u\in {\cal U}$ we have for $x\in E$ that $\|P_{2^{-m}}^{(m),u_n(x)}(x,\cdot)-P_{2^{-m}}^{(m),u(x)}(x,\cdot)\|_V\to 0$
as $n\to \infty$.}

\medskip

By analogy to Proposition 1 and also Proposition 2 of \cite{LS0} we have by induction

\noindent
{\bf Lemma 3.} {\it Under (uCont) for $u_n\to u\in {\cal U}$ and any $k\in N$ we have
\begin{equation}
\|P_{k2^{-m}}^{(m),u_n(x)}(x,\cdot)-P_{k2^{-m}}^{(m),u(x)}(x,\cdot)\|_V\to 0.
\end{equation}
as $n\to \infty$.} 

Our main result can be formulated as follows

\noindent
{\bf Theorem 2.} {\it Under (uCont), (UEd) and (FPV) we have that 
\begin{equation}\label{invconv2}
\|\mu_m^{u_n}-\mu_m^u\|_V \to 0
\end{equation}
as $n\to \infty$. Additionally under (CCon), (ER) and (ERd) we have that 
\begin{equation}\label{confunc}
J_x^m(u_n)\to J_x(u)
\end{equation}
as $n, m\to \infty$. Moreover
\begin{equation}\label{confunc1}
J_x(u_n)\to J_x(u)
\end{equation}
as $n\to \infty$.}
  
{\bf Proof.} To prove (\ref{invconv2}) we use  Lemma 1 and then Lemma 3. Convergence (\ref{confunc}) follows from 
 (\ref{invconv2}) and Theorem 1. Convergence (\ref{confunc1}) can be shown from Lemma 3, Lemma 1 and Corollary 2. 
 
\rightline{$\Box$}

\section{Risk sensitive control}

We shall assume that 

\medskip 

\noindent
(uUE) {\it for each $u\in {\cal U}$ there is $\Delta_u\in (0,1)$ such that we have  $\sup_{m\in N} \sup_{x, x'\in E} \sup_{B\in {\cal E}} P_1^{(m),u}(x,B)-P_1^{(m),u}(x',B):=\Delta_u<1$.}

\medskip
\noindent
It is clear that under (uUE) Markov process $(X_n^{(m),u})$ has a unique invariant measure $\mu_m^u$ (see \cite{Doob}). Furthermore additionally under (Conv) with $V\equiv 1$ we have that 
\begin{equation}\label{cuerg}
\sup_{x, x'\in E} \sup_{B\in {\cal E}} P_1^{u}(x,B)-P_1^{u}(x',B)\leq\Delta_u<1.
\end{equation}
Then process  $(X_n^{u})$ has a unique invariant measure $\mu^u$. 

We also assume that 

\medskip

\noindent
(uEquiv) {\it for each $u\in {\cal U}$ there is $k\in N$ such that we have that $\sup_{m\in N} \sup_{x, x'\in E} \sup_{B\in {\cal E}} 
{P_k^{(m),u}(x,B) \over P_k^{m(m),u}(x',B)}:=K_u<\infty$}. 

\medskip

\noindent
Under (Conv) and (uEquiv) we have that 
\begin{equation}\label{uEqu}
\sup_{x, x'\in E} \sup_{B\in {\cal E}} 
{P_k^{u}(x,B) \over P_k^{u}(x',B)}\leq K_u<\infty.
\end{equation}
  
\noindent
{\bf Example 2.} Assume that diffusion process $(X_t^u)$ defined in Example 1 is reflected in a bounded regular domain. Then following Theorem 2.1 of \cite{MR} (see also \cite{GM}) we can show property (\ref{cuerg}). Since transition densities are bounded away from zero we also have that (\ref{uEqu}) is satisfied.

Let $B(E)$ be the set of bounded Borel measurable functions on $E$ with supremum norm. For $g\in B(E)$ define so called span norm  $\|g\|_{sp}=\sup_{x\in E} g(x)- \inf_{x'\in E} g(x')$. For $u\in {\cal U}$ and $f, g\in B(E)$, and $\alpha\in (-\infty, + \infty)\setminus \left\{0\right\}$ define

\begin{equation}
\Psi^{(m), u,\alpha} g(x)= {1 \over \alpha}\ln E_x^u\left\{\exp\left\{\alpha \sum_{i=0}^{2^m-1}2^{-m}c(X_{i 2^{-m}}^{(m),u},u( X_{i 2^{-m}}^{(m),u})) +\alpha g(X_{1}^{(m),u})\right\}\right\}. 
\end{equation}   

We have

\noindent
{\bf Theorem 3.} {\it  Under (uUE) for $\alpha\neq 0$ the operator $\Psi^{(m), u,\alpha}$ is a local contraction in the span norm in the space $B(E)$ for $u\in {\cal U}$, i.e. there is a function $\gamma_\alpha: (0,\infty) \mapsto [0,1)$, which does not depend on $m$, such that whenever for $g_1,g_2 \in B(E)$ we have $\|g_1\|_{sp}\leq M$ and $\|g_2\|_{sp}\leq M$ then
\begin{equation}\label{lcont}
\|\Psi^{(m), u,\alpha} g_1 - \Psi^{(m), u,\alpha} g_2\|_{sp}\leq \gamma_\alpha(M)\|g_1-g_2\|_{sp}.
\end{equation}
Furthermore additionally under (uEquiv) the $k$-th iteration of  $\Psi^{(m),u,\alpha}$ transforms the space $B(E)$ to the subspace of $B(E)$ with the span norm less than $\tilde{K}_u$, with $\tilde{K}_u$ depending on $K_u$ from (uUE). Consequently $\Psi^{(m),u,\alpha}$ after $k$-th iteration is a global contraction.} 

\noindent
{\bf Proof.}  Local contractivity follows from Theorem 3, Corollary 4 and 5 in \cite{LS2} in a similar way as in section 2 of  \cite{LS4}. We give here only few hints. Using dual representation of the operator $\Psi$ (see Proposition 1.42 of \cite{DE}) we have that for $\alpha<0$

\begin{eqnarray}\label{er1}
&&\Psi^{(m), u,\alpha} g(x)= inf_{\nu\in P_x(D_E[0,1])}  \int_{D_E[0,1]} \left(2^{-m}\sum_{i=0}^{2^m-1}c(z_{i 2^{-m}},u(z_{i 2^{-m}}))\right. + \nonumber \\
&&\left. g(z_1)\right)\nu(dz)-{1\over \alpha} H(\nu, P_{[0,1]}^{(m),u}(x,\cdot)) 
\end{eqnarray}
and for $\alpha>0$ 
\begin{eqnarray}\label{er2}
&&\Psi^{(m), u,\alpha} g(x)= sup_{\nu\in P_x(D_E[0,1])}  \int_{D_E[0,1]} \left(2^{-m}\sum_{i=0}^{2^m-1}c(z_{i 2^{-m}},u(z_{i 2^{-m}}))\right. + \nonumber \\
&&\left. g(z_1)\right)\nu(dz)-{1\over \alpha} H(\nu, P_{[0,1]}^{(m),u}(x,\cdot)), 
\end{eqnarray}
where $D_E[0,1]$ is the set of all c\`adl\`ag trajectories on the time interval $[0,1]$, while $P_x(D_E[0,1])$ is the set of probability measures on the set $D_E[0,1]$ starting from $x\in E$ and $H$ denotes entropy between measures $\nu$ and  $P_{[0,1]}^{(m),u}(x,\cdot)$ defined as follows $H(\nu_1,\nu_2):=\int_{D_E[0,1]} \ln({d\nu_1 \over d \nu_2})d\nu_1$ when $\nu_1$  
is absolutely continuous with respect to $\nu_2$, and is equal to $+\infty$ otherwise.  Infimum in (\ref{er1}) or supremum in (\ref{er2}) is attained by the measure on $D_E[0,1]$ of the form

\begin{eqnarray}
&&\nu_{x,\alpha g}^{(m),u}(dz):=\exp\left(\alpha \sum_{i=0}^{2^m-1}2^{-m}c(z_{i 2^{-m}},u(z_{i 2^{-m}}))   + \alpha g(z_1)\right)P_{[0,1]}^{(m),u}(x,dz) \nonumber \\
&&\left[E_x^u\left\{\exp\left\{\alpha \sum_{i=0}^{2^m-1}2^{-m}c(X_{i 2^{-m}}^{(m),u},u( X_{i 2^{-m}}^{(m),u})) +\alpha g(X_{1}^{(m),u})\right\}\right\} \right]^{-1}.
\end{eqnarray}
Define now the measure on ${\cal E}$ 
\begin{eqnarray}   
&&\overline{\nu_{x,\alpha g}^{(m),u}}(B):= \left[E_x^u\left\{1_B(X_1^{(m),u}) \right. \right. \exp\left\{\alpha \sum_{i=0}^{2^m-1}2^{-m}c(X_{i 2^{-m}}^{(m),u},u( X_{i 2^{-m}}^{(m),u})) \right. \nonumber \\
&& \left.\left. \left.  +\alpha g(X_{1}^{(m),u})\right\}\right\}\right] \left[E_x^u\left\{\exp\left\{\alpha \sum_{i=0}^{2^m-1}2^{-m}c(X_{i 2^{-m}}^{(m),u},u( X_{i 2^{-m}}^{(m),u}))\right.\right. \right. \nonumber \\
&& \left.\left. \left.  +\alpha g(X_{1}^{(m),u})\right\}\right\} \right]^{-1}.
\end{eqnarray}

For $g_1,g_2\in B(E)$ and $x_1,x_2\in E$ and $\alpha<0$ using (\ref{er1})-(\ref{er2}) we obtain
\begin{eqnarray}\label{ineqimp}
&& \Psi^{(m),u,\alpha}g_1(x_1)-\Psi^{(m), u,\alpha}g_2(x_1)-\Psi^{(m),u,\alpha}g_1(x_2)+ \nonumber \\
&&\Psi^{(m), u,\alpha}g_1(x_2) \leq\|g_1-g_2\|_{sp} \sup_{B\in {\cal E}} (\nu_1-\nu_2)(B)
\end{eqnarray}
where $\nu_1:=\overline{\nu_{x_1,\alpha g_2}^{(m),u}}$ and 
$\nu_2:=\overline{\nu_{x_2,\alpha g_1}^{(m),u}}$. In the case of $\alpha>0$ we replace $g_1$ with $g_2$ in definitions of $\nu_1$ and $\nu_2$. Assume now that for $x_n, x_n'\in E$, $g_{1,n},g_{2,n}$, such that $\|g_{1,n}\|_{sp}\leq M$ and $\|g_{2,n}\|_{sp}\leq M$ in the place of $x_1,x_2$, $g_1,g_2$ and $m_n\in N$ in the place of $m$, $B_n\in {\cal E}$ we have $\nu_1(B_n)\to 1$ and $\nu_2(B_n)\to 0$.
Then in the case of $\alpha<0$ we obtain $\nu_2(B_n)\geq P_1^{(m_n),u}(x_{2,n},B_n)e^{-\|c\|_{sp}-M}$ 
and $\nu_1(B_n^c)\geq P_1^{(m_n),u}(x_{1,n},B_n^c)e^{-\|c\|_{sp}-M}$, 
which implies that $P_1^{(m_n),u}(x_{2,n},B_n)\to 0$ and 
$ P_1^{(m_n),u}(x_{1,n},B_n^c)\to 0$, as $n\to \infty$ contradicting 
(uUE). In the case of $\alpha>0$ we have a similar contradiction. 
This completes the proof of local contractivity of $\Psi^{(m), u,\alpha}$     
with a  Lipschitz constant $\gamma_\alpha(M)$.  
Global contraction then follows from Remark 4 and Proposition 6 in \cite{LS1}. Namely, under (uEquiv) we have
\begin{equation}\label{bound1}
\|(\Psi^{(m),u,\alpha})^k g\|_{sp} \leq k \|c\|_{sp} + \ln K_u.
\end{equation}
This means that $k$-th iteration of the operator $\Psi^{(m), u,\alpha}g$ no matter what $g\in B(E)$ was chosen has a uniformly bounded span norm. 
The proof of Theorem 3 is therefore completed.

\rightline{$\Box$}

Basing on Theorem 3 we obtain the solutions to certain versions of the Poisson equations

\noindent
{\bf Corollary 3.} {\it Under assumptions of Theorem 3  for $u\in {\cal U}$ there is a constant $\lambda^{(m), u,\alpha}$ and a function $w^{(m), u,\alpha}\in B(E)$ such that for $x\in E$ we have
\begin{eqnarray}\label{riskbel}
&&e^{\alpha w^{(m),u,\alpha}(x)} =  E_x^u\left\{\exp\left\{\alpha \sum_{i=0}^{2^m-1}2^{-m}(c(X_{i 2^{-m}}^{(m),u},u(X_{i 2^{-m}}^{(m),u}))-  \right. \right. \nonumber \\
&&\left. \left. \lambda^{(m),u,\alpha})+\alpha w^{(m),u,\alpha}(X_{1}^{(m),u})\right\}\right\} 
\end{eqnarray}
Moreover $\|w^{(m), u,\alpha}\|_{sp}\leq \tilde{K}_u$, where $\tilde{K}_u$ depends on  $K_u$ from (uEquiv) and the function $\gamma_\alpha$.}

\noindent
{\bf Proof.} By Theorem 1 there is a fixed point $w^{(m),u,\alpha}$ of the operator $\Psi^{(m),u,\alpha}$ i.e. $\|\Psi^{(m),u,\alpha} w^{(m),u,\alpha}-w^{(m),u,\alpha}\|_{sp}=0$. Therefore there is a constant $\lambda^{(m),u,\alpha}$ such that $\Psi^{(m),u,\alpha} w^{(m),u,\alpha}(x)-\lambda^{(m),u,\alpha} = w^{(m),u,\alpha}$, which completes the proof.

\rightline{$\Box$} 

\noindent
{\bf Corollary 4.}\label{Cor4} {\it If $\lambda^{(m), u,\alpha}$ and a function $w^{(m), u,\alpha}\in B(E)$ such that  $\|w^{(m), u,\alpha}\|_{sp}\leq \tilde{K}_u$ are solutions to the equation (\ref{riskbel}) then for any $k\in N$ we have
\begin{equation}\label{driskbel}
e^{\alpha w^{(m),u,\alpha}(x)} =  E_x^u\left\{\exp\left\{\alpha \sum_{i=0}^{k2^m-1}2^{-m}(c(X_{i 2^{-m}}^{(m),u},u(X_{i 2^{-m}}^{(m),u}))-  \lambda^{(m),u,\alpha})+\alpha w^{(m),u,\alpha}(X_{k}^{(m),u})\right\}\right\} 
\end{equation}
and consequently 
\begin{equation}\label{ddriskbel}
| \lambda^{(m),u,\alpha}-{1\over \alpha}{1\over k}\ln E_x^u\left\{\exp\left\{\alpha \sum_{i=0}^{k2^m-1}2^{-m}(c(X_{i 2^{-m}}^{(m),u},u(X_{i 2^{-m}}^{(m),u}))\right\}\right\}|\leq 2{\tilde{K}_u \over k}.
\end{equation}
Therefore for any $x\in E$ we have
\begin{equation}\label{driskopt}
\lambda^{(m),u,\alpha}=I_x^{\alpha,2^{-m}}(u).
\end{equation}}

\noindent{\bf Proof.} Iterating equation (\ref{riskbel}) we obtain (\ref{driskbel}). Taking into account that  $\|w^{(m), u,\alpha}\|_{sp}\leq \tilde{K}_u$ we then obtain (\ref{ddriskbel}). Since $c$ is bounded we have that
\begin{eqnarray}
&&\liminf_{k\to \infty} {1\over \alpha}{1\over k 2^m} \ln E_x^u\left\{\exp\left\{\alpha 2^{-m}\sum_{i=0}^{k-1}(c(X_{i 2^{-m}}^{(m),u},u(X_{i 2^{-m}}^{(m),u}))\right\}\right\}=\nonumber \\
&&\liminf_{k\to \infty} {1\over \alpha}{1\over k} \ln E_x^u\left\{\exp\left\{\alpha 2^{-m}\sum_{i=0}^{k2^m-1}(c(X_{i 2^{-m}}^{(m),u},u(X_{i 2^{-m}}^{(m),u}))\right\}\right\},
\end{eqnarray}
from which (\ref{driskopt}) follows. 

\rightline{$\Box$} 

Assume now

\medskip

\noindent
(eConv) {\it for each $u\in {\cal U}$ measures 
\begin{equation*}
{\cal M}_x^{(m),u,\alpha}(B):=E_x^u\left\{\exp\left\{\alpha \sum_{i=0}^{2^m-1}2^{-m}(c(X_{i 2^{-m}}^{(m),u},u(X_{i 2^{-m}}^{(m),u})) \right\} 1_B(X_{1}^{(m),u}) \right\}
\end{equation*}
defined for $B\in {\cal E}$ converge in variation norm to the measure 
\begin{equation*}
{\cal M}_x^{u,\alpha}(B):=E_x^u\left\{\exp\left\{\alpha \int_0^1 (c(X_{s}^{u},u(X_{s}^{u}))ds\right\} 1_B(X_{1}^{u})\right\},
\end{equation*}
 as $m\to \infty$. }
\medskip

We then have

\noindent
{\bf Theorem 4.} {\it Under (eConv),  (uUE) and (uEquiv) for $u\in {\cal U}$ and there is a constant $\lambda^{u,\alpha}$ and a function $w^{u,\alpha}\in B(E)$ such that $\|w^{u,\alpha}\|_{sp}\leq \tilde{K}$ and for $x\in E$ we have
\begin{equation}\label{criskbel}
e^{\alpha w^{u,\alpha}(x)} =  E_x^u\left\{\exp\left\{\alpha \int_0^1 (c(X_{s}^{u},u(X_{s}^{u}))ds-\lambda^{u,\alpha})+\alpha w^{u,\alpha}(X_{1}^{u})\right\}\right\}. 
\end{equation}}

\noindent
{\bf Proof.} 
Let for $g\in B(E)$
\begin{equation}
\Psi^{u,\alpha} g(x)=  E_x^u\left\{\exp\left\{\alpha \int_0^1 (c(X_{s}^{u},u(X_{s}^{u}))ds+\alpha w^{u,\alpha}(X_{1}^{u})\right\}\right\}.\nonumber
\end{equation}

By (\ref{ineqimp}) since $\Psi^{(m),u,\alpha}g_i(x_j) \to \Psi^{u,\alpha}g_i(x_j)$ for $i,j\in \left\{1,2\right\}$, as $m\to \infty$, we obtain 
\begin{equation}
\|\Psi^{u,\alpha} g_1 - \Psi^{u,\alpha} g_2\|_{sp}\leq \gamma_\alpha(M)\|g_1-g_2\|_{sp},
\end{equation} 
for $\|g_1\|_{sp}\leq M$ and $\|g_2\|_{sp}\leq M$ with the same $\gamma_\alpha(M)$ as in (\ref{lcont}).
Letting now $m\to \infty$ in (\ref{bound1}), taking into account that a version of (eConv) also holds for time $k$ (instead of $1$) we obtain the bound for iterations of $\Psi^{u,\alpha} g$ with $g\in B(E)$, from which existence of a unique fixed point of  $\Psi^{u,\alpha}$ with suitable bound follows. 

\rightline{$\Box$} 

In analogy to Corollary 4 we obtain

\noindent
{\bf Corollary 5.}\label{Cor5} {\it If $\lambda^{u,\alpha}$ and a function $w^{u,\alpha}\in B(E)$ such that  $\|w^{ u,\alpha}\|_{sp}\leq \tilde{K}_u$ are solutions to the equation (\ref{criskbel}) then for any $k\in N$ we have
\begin{equation}\label{cdriskbel}
e^{\alpha w^{u,\alpha}(x)} = E_x^u\left\{\exp\left\{\alpha (\int_0^k (c(X_{s}^{u},u(X_{s}^{u}))ds- \lambda^{u,\alpha}) +\alpha w^{u,\alpha}(X_{k}^{u})\right\}\right\} 
\end{equation}
and consequently 
\begin{equation}\label{cddriskbel}
| \lambda^{u,\alpha}-{1\over \alpha}{1\over k}\ln E_x^u\left\{\exp\left\{\alpha \int_0^k c(X_{s}^{u},u(X_{s}^{u}))ds \right\} \right\}| \leq 2{\tilde{K}_u\over k}.
\end{equation}
Therefore for any $x\in E$ we have
\begin{equation}\label{cdriskopt}
\lambda^{u,\alpha}=I_x^{\alpha}(u).
\end{equation}}

\noindent
{\bf Proof.} Note that  (\ref{cdriskbel}) and (\ref{cddriskbel}) follow easily from (\ref{criskbel}) (similarly as in the proof of Corollary 4). To show (\ref{cdriskopt}) it sufficies to notice that by boundedness of $c$ and  (\ref{cddriskbel}) we have
\begin{equation}
 \liminf_{t\to \infty} {1\over \alpha} {1\over t} \ln E_x^u\left\{e^{\alpha \int_0^t c(X_s^u,u(X_s^u))ds}\right\}= \nonumber \\
\liminf_{k\to \infty} {1\over \alpha} {1\over k} \ln E_x^u\left\{e^{\alpha \int_0^k c(X_s^u,u(X_s^u))ds}\right\},
\end{equation}
where first line we have a limit of positive real $t$ going to $\infty$, while in the second line over positive integer $k$ going to $\infty$. 

\rightline{$\Box$}
   
The following Corollary summarizes just obtained results

\noindent
{\bf Corollary 6.} 
{\it Under (eConv),  (uUE) and (uEquiv)  $u\in {\cal U}$ we have
\begin{equation}
I_x^{\alpha, 2^{-m}}(u)\to I_x^{\alpha}(u),
\end{equation}
as $m\to \infty$.}

\noindent
{\bf Proof.} Clearly by Corollaries 4 and 5 we have that $I_x^{\alpha, 2^{-m}}(u)=\lambda^{(m),u,\alpha}$ and $I_x^{\alpha}(u)=\lambda^{u,\alpha}$. Now using (eConv) to (\ref{driskbel}) and (\ref{ddriskbel}) we obtain that  $\lambda^{(m),u,\alpha}\to \lambda^{u,\alpha}$, as $m\to \infty$.

\rightline{$\Box$}

\noindent
{\bf Remark 1.} {\it Assumption (uUE) plays an important role to study discrete time risk sensitive Bellman equation.
We require it to be satisfied uniformly with respect to discretization step, which is important when we let discretization step converging to $0$.  Discrete time risk sensitive problems can be also studied using splitting technics as in the paper \cite{DiMS2}. This however would require a number of additional assumptions. 
Assumption (uEquiv) can be replaced by requiring small risk $|\alpha|$ as was studied in the papers \cite{DiMS3} or \cite{PS}. 
Using assumption (uUE) we are looking for a bounded solution to (\ref{riskbel}) (see \cite{DiMS1}). We can use also other technics based on Krein Rutman theorem (see \cite{LS3} and \cite{A}) or suitable Lyapunov conditions (see \cite{Ch}) and work with unbounded solutions. In such case we shall also require more assumptions.} 

We now consider stability of functional $I_x^u$. We have

\noindent{\bf Theorem 5.} {\it Assume (eConv), (uUE), (uEquiv) are satisfied for each $u\in {\cal U}$ with $\sup_{u\in {\cal U}} \Delta_u<1$, $\sup_{u\in {\cal U}}K_u<\infty$. Then under (uCont) for ${\cal U}\ni u_n\to u\in {\cal U}$ as $n\to \infty$ we have for each $m\in N$ 
\begin{equation}
 I_x^{\alpha, 2^{-m}}(u_n)=\lambda^{(m), u_n,\alpha}\to  I_x^{\alpha, 2^{-m}}(u)=\lambda^{(m),u,\alpha}.
\end{equation}
Furthermore, when additionally (eConv) is satisfied uniformly for $(u_n)$ we have
\begin{equation} \label{finconv}
 I_x^\alpha(u_n)=\lambda^{u_n,\alpha}\to  I_x^\alpha(u)=\lambda^{u,\alpha}, 
\end{equation}
as $n\to \infty$.} 

\noindent{\bf Proof.} We have 
\begin{equation}
|\lambda^{u_n,\alpha}-\lambda^{u,\alpha}|\leq |\lambda^{u_n,\alpha}-\lambda^{(m),u_n,\alpha}|+ |\lambda^{(m),u_n,\alpha}-\lambda^{(m),u,\alpha}|+|\lambda^{(m),u,\alpha}-\lambda^{u,\alpha}|.
\end{equation}
Now by (\ref{ddriskbel}) we obtain
\begin{equation}
|\lambda^{(m),u_n,\alpha}-\lambda^{(m),u,\alpha}|\leq 4 \sup_{u\in {\cal U}}{\tilde{K}_u\over k} + {1\over \alpha}{1\over k}
W((m),u_n,u,\alpha,k),
\end{equation}
where
\begin{eqnarray}
&&W((m),u_n,u,\alpha,k):= |\ln E_x^{u_n}\left\{\exp\left\{\alpha \sum_{i=0}^{k2^m-1}2^{-m} 
(c(X_{i 2^{-m}}^{(m),u_n}, u(X_{i 2^{-m}}^{(m),u_n}))\right\}\right\}- \nonumber \\
&&\ln E_x^{u}\left\{\exp\left\{\alpha \sum_{i=0}^{k2^m-1}2^{-m}  
(c(X_{i 2^{-m}}^{(m),u},u(X_{i 2^{-m}}^{(m),u}))\right\}\right\}|. 
 \end{eqnarray}
It is clear that under (uCont) for each $m\in N$ and $k\in N$, $W((m),u_n,u,\alpha,k)$ converges to $0$ as $n\to \infty$. Furthermore $\sup_{u\in {\cal U}} \tilde{K}_u<\infty$. Consequently letting first $n\to \infty$ then $k\to \infty$ we obtain that $|\lambda^{(m),u_n,\alpha}-\lambda^{(m),u,\alpha}|\to 0$ as $n\to \infty$. Using Corollary 6 and the fact that (eConv) is satisfied uniformly for $(u_n)$ we obtain (\ref{finconv}). 
 
\rightline{$\Box$}
    
\section{CONCLUSIONS}
In the paper we justify the use of natural approximation procedure for continuous time controlled Markov processes over long time horizon. Namely, instead of using Markov control $u(X_t)$ at each time $t$ we choose control $u(X_{nh})$ at times $nh$ and consider control fixed  in the time intervals $[nh,(n+1)h)$.  It appears that under reasonable assumptions we obtain a good approximation of the average reward per unit time functional as well as long run risk sensitive functional.  This way we obtain a feasible construction of nearly optimal controls for continuous time controlled Markov processes, which can be used in various applications.

\end{document}